\documentclass[12pt,leqno]{article}
\usepackage[pdftex]{graphicx}
\usepackage{amsfonts,amssymb,amsmath}
\usepackage{ccaption,mathrsfs,amsmath,amssymb,amsthm}
\usepackage{amsmath, amscd, amsfonts,  tikz, times}

\def\tcr{\textcolor{red}}
\def\tcb{\textcolor{blue}}
\newcounter{conjecture}
\newcounter{remark}

\newcommand{\eqnsection}{
\renewcommand{\theequation}{\thesection.\arabic{equation}}
\makeatletter
\csname @addtoreset\endcsname{equation}{section}
\makeatother}
\newtheorem{theorem}{Theorem}

\newtheorem{lemma}{Lemma}
\newtheorem{defn}{Definition}

\newtheorem{cor}{Corollary}

\newcommand{\dd}{\delta}

\newcommand{\lar}{\longrightarrow}

\newcommand{\aaa}{\alpha}

\newcommand{\CC}{\mathbb{C}}
\newcommand{\ZZ}{\mathbb{Z}}

\def \ov{\overline}

\def \be{\begin{equation}}
\def \ee{\end{equation}}
\def \bt{\begin{theorem}}
\def \et{\end{theorem}}
\def \bea{\begin{eqnarray}}
\def \eea{\end{eqnarray}}
\def \bas{\begin{eqnarray*}}
\def \eas{\end{eqnarray*}}

\newcommand {\rrr}[1]{(\ref{#1})}

\def \Om{\Omega}

\def \si{\sigma}

\def \th{\theta}

\def \ff{\infty}

\def \AA{{\cal A}}

\def \DD{{\mathbb D}}

\def \GG{{\cal G}}
\def \HH{{\cal H}}
\def \HHH{{\mathbb H}}

\def \RR{{\mathbb R}}

\def \({\left(}
\def \){\right)}

\newcommand{\fr}{\frac}

\def \vski{\vspace{12pt}}

\newcommand{\bsh}{\backslash}

\def \bc{\begin{center} }
\def \ec{\end{center} }
\def \bs{\begin{slide} }
\def \es{\end{slide} }

\def\square{{\vcenter{\vbox{\hrule height.3pt
         \hbox{\vrule width.3pt height5pt \kern5pt
            \vrule width.3pt}
         \hrule height.3pt}}}}
\def\qed{{\hfill $\Box$ \bigskip}}

\newcounter{cccases}
\eqnsection

\begin{document}
\title{Complex analytic proofs of two probabilistic theorems}
\author{\textit{Greg Markowsky and Clayton McDonald} 
\thanks{Email: greg.markowsky@monash.edu, clayton.mcdonald@monash.edu }\\
\textit{School of Mathematics, Monash University}
}

\bibliographystyle{amsplain}
\maketitle \eqnsection \setlength{\unitlength}{2mm}

\begin{abstract}
\noindent In this paper, we use purely complex analytic techniques to prove two results of the first author which were hitherto given only probabilistic proofs. 

A general form of the Phragm\'en-Lindel\"of principle states that if the $p$\textsuperscript{th} Hardy norm of the conformal map from the disk to a simply connected domain is finite, then an analytic function on that domain is either bounded by its supremum on the boundary or else goes to $\ff$ along some sequence more rapidly than $e^{|z|^{p}}$. We will prove this and discuss a number of special cases.

We also derive a series expansion for the Green's function of a disk, and show how it leads to an infinite product identity. The celebrated infinite product expansions for sine and cosine are realized as special cases.

\noindent

\vski


\vski

\end{abstract}

\section{Introduction}

In two past papers by the first author \cite{greggreenstoppingtimes, mebabyme}, several complex analytic results were proved using the conformal invariance of Brownian motion; however, they could be stated in equivalent form without any reference to probablility. The object of this paper is to show how they can be proved without Brownian motion, using only analysis.

\section{The Phragm\'en-Lindel\"of principle}

The Phragm\'en-Lindel\"of principle is a method by which the maximum modulus principle can be generalized to certain unbounded domains in $\CC$. The principle roughly states that, on particular domains, analytic functions must either be bounded by their supremum on the boundary of the domain or tend rapidly to $\ff$ along some sequence. In practice, when the principle is discussed precise statements are often made and proved for specific domains, such as the angular wedge $N_\aaa = \{re^{i\th}: r \in (0,\ff), \th \in (\frac{-\aaa}{2},\frac{\aaa}{2})\}$ (see for instance \cite[Sec. VI.4]{conway}). However, in \cite{mebabyme} a general result was proved using probabilistic arguments which applies to a number of special cases, including $N_\aaa$. Although the proof given was quite short, it depended upon familiarity with the properties of Brownian motion, and may therefore be somewhat inaccessible for analysts without a background in stochastic processes. A purely analytic proof is therefore desirable, and the purpose of this note is to provide one.\\

Although the theorem proved in \cite{mebabyme} deals with more general domains, we will restrict ourselves to simply connected ones in this paper. If $W \subsetneq \CC$ is simply connected then the Riemann Mapping Theorem guarantees the existence of a conformal map $\phi_a$ mapping $\DD$ onto $W$ which takes $0$ to $a$. The Hardy norm $|| \cdot ||_{H^{p}}$ of $\phi_a$ is then defined as

\begin{equation} \label{h2def}
||\phi_a||_{H^{p}} := \Big(\sup_{r<1} \frac{1}{2\pi} \int_{0}^{2\pi} |\phi_a(re^{i \th})|^{p} d\th \Big)^{1/p} .
\end{equation}

The map $\phi_a$ is not uniquely determined; however any two such maps differ only by precomposition with a rotation, so the value of $||\phi_a||_{H^{p}}$ is independent of the choice of $\phi_a$. We set $\HH^a_{p}(W) = ||\phi_a||_{H^{p}}$. It may be shown that if $a,b \in W$ then $\HH^a_{p}(W), \HH^b_{p}(W)$ are either both finite or both infinite. In light of this, let us write $\HH_{p}(W) < \ff$ whenever $\HH^a_{p}(W) < \ff$ for some (and thus every) $a \in W$. Our result is as follows.

\begin{theorem} \label{phraglind1}
Let $W$ be a simply connected domain with $\HH_{p}(W) < \ff$. Suppose that $f$ is an analytic function on $W$ such that $\limsup_{z \lar \dd W}|f(z)| \leq K < \ff$, and $|f(z)| \leq Ce^{C|z|^{p}}$ for some $C>0$. Then $|f(z)| \leq K$ for all $z \in W$.
\end{theorem}

Before giving the proof in Section \ref{proof}, for completeness we briefly discuss the special cases listed in \cite{mebabyme}. To obtain our first variant of Theorem \ref{phraglind1}, we remark that it is known that if $W \subsetneq \CC$ is simply connected then $\HH_{p}(W)<\ff$ for any $p < \frac{1}{2}$. This is proved in \cite{burk}, and follows easily as well from the main result of \cite{burn}. Furthermore, the fact that the $H^{p}$-Hardy norm of the Koebe function $f(z)=\frac{z}{(1-z)^2}$, which maps $\DD$ conformally onto $\CC \bsh (-\ff, 1/4]$, is finite if and only if $p<1/2$ shows it cannot be improved. We obtain

\begin{cor} \label{simpcon}
Suppose that $f$ is an analytic function on a simply connected domain $W \subsetneq \CC$ such that $\limsup_{z \lar \dd W}|f(z)| \leq K < \ff$, and $|f(z)| \leq Ce^{C|z|^{p}}$ for some $p < \frac{1}{2}$. Then $|f(z)| \leq K$ for all $z \in W$.
\end{cor}

A domain $W$ is {\it spiral-like of order $\si \geq 0$ with center $a$} if, for any $z \in W$, the spiral $\{a+(z-a) \mbox{ exp}(te^{-i \si}) : t \leq 0\}$ also lies within $W$ (cf. \cite{space}). In \cite{hansenspi}, Hansen gave a geometric condition for the finiteness of $\HH_{p}(W)$ for any spiral-like domain $W$, which we now describe. There is no loss of generality in assuming $a=0$, and we will do so henceforth. Hansen showed that the key quantity for our purposes is the measure of the largest arc in the set $W \cap \{|z|=r\}$ (taken as a set on the circle). Set

\begin{equation} \label{bigmax}
\AA_{r,W} = \max \{m(E): E \mbox{ is a subarc of } W \cap \{|z|=r\}\},
\end{equation}

where $m$ denotes angular Lebesgue measure on the circle. Spiral-likeness implies that $\AA_{r,W}$ is nondecreasing in $r$, so we may let $\AA_W = \lim_{r \nearrow \ff} \AA_{r.W}$. Hansen's result is as follows.

\begin{theorem} \label{lamoreau}
If $W$ is a spiral-like domain of order $\si$ with center 0, then $\HH_{p}(W) < \ff$ if, and only if, $p < \frac{\pi}{\AA_W \cos^2 \si}$.
\end{theorem}

\vspace{.02in}

\begin{cor} \label{page1}
Suppose $W$ is a spiral-like domain of order $\si$ with center $0$. If $f$ is an analytic function on $W$ such that $\limsup_{z \lar \dd W}|f(z)| \leq K < \ff$, and $|f(z)| \leq Ce^{C|z|^{p}}$ for any $p < \frac{\pi}{\AA_W \cos^2 \si}$, then $|f(z)| \leq K$ for all $z \in W$.
\end{cor}

A domain $W$ is called {\it star-like with center $a$} if the line segment connecting $a$ to $z$ lies within $W$ for every $z \in W$. In addition to being of interest in their own right, star-like domains often figure prominently in a student's first exposure to complex analysis, as a number of fundamental results (such as Cauchy's Theorem) are commonly proved first for star-like domains before being extended to more general ones. Note that a star-like domain is simply spiral-like of order $\si =0$, and Corollary \ref{page1} therefore takes the following form as a special case.

\begin{cor} \label{kristen}
Suppose $W$ is a star-like domain with center $0$. If $f$ is an analytic function on $W$ such that $\limsup_{z \lar \dd W}|f(z)| \leq K < \ff$, and $|f(z)| \leq Ce^{C|z|^{p}}$ for any $p < \frac{\pi}{\AA_W}$, then $|f(z)| \leq K$ for all $z \in W$.
\end{cor}

Note that convex domains are trivially star-like, and the previous theorem therefore applies to any convex domains as well. Let us now set $N_\aaa = \{re^{i\th}: r \in (0,\ff), \th \in (\frac{-\aaa}{2},\frac{\aaa}{2})\}$; $N_\aaa$ is the angular wedge with vertex at 0 which is symmetric about the real axis and has angular width $\aaa$. $N_\aaa$ is star-like, and Corollary \ref{kristen} therefore reduces further to the following.

\begin{cor} \label{edgewedge}
Suppose that $f$ is an analytic function on $N_\aaa$ such that $\limsup_{z \lar \dd W}|f(z)| \leq K < \ff$, and $|f(z)| \leq Ce^{C|z|^{p}}$ for some $p < \frac{\pi}{\aaa}$. Then $|f(z)| \leq K$ for all $z \in N_\aaa$.
\end{cor}

This is probably the most commonly stated form of the Phragm\'en-Lindel\"of principle, and appears in many different texts (such as \cite[Sec. VI.4]{conway}). It may also be noted that this example shows that Corollary \ref{edgewedge} (as well as Corollaries \ref{kristen} and \ref{page1}, and, in turn, Theorem \ref{phraglind1}) is sharp in the sense that the function $f(z) = e^{z^{\pi/\aaa}}$ is analytic on $N_\aaa$ and bounded in modulus by 1 on $\dd N_\aaa$, but is clearly unbounded on $N_\aaa$.\\

\section{Proof of Theorem \ref{phraglind1}} \label{proof}

The key is the following lemma.

\begin{lemma} \label{calig}
Let $W$ be a domain with $\HH_{p}(W) < \ff$. Suppose that $u$ is a continuous function on $cl(W)$ which is subharmonic on $W$ and satisfies $\sup_{z \in \dd W} u(z) \leq K$, for some $K > 0$. Suppose further that $u(z) \leq C|z|^{p} + C$ for some $C < \ff$. Then $u(z) \leq K$ for all $z \in W$.
\end{lemma}

{\bf Proof:} Fix $a \in W$, and let $\phi_a: \DD \lar W$ be as above. Since $u$ is subharmonic, the composition $u \circ \phi_a$ is as well, and we therefore have for any $r<1$

\begin{equation} \label{nowayout}
u(a) = u(f(0)) \leq \frac{1}{2\pi} \int_{0}^{2\pi} u(\phi_a(re^{i\th}))d\th.
\end{equation}

Now we make use of the fact that $\HH_{p}(W) < \ff$. Standard results in the theory of Hardy spaces show that the functions $\phi_a(re^{i\th})$ approach a limit function $\phi_a(e^{i\th})$ in $L^p([0,2\pi])$ for almost every $\th$ in $[0,2\pi]$ as $r \nearrow 1$, and $\phi_a(e^{i\th})$ takes values in $\dd W$ a.e. (see \cite[Thm. 17.10]{rud}). Furthermore, an inequality due to Hardy and Littlewood (see \cite{hardwood}) allows us to conclude that the function $M_{\phi_a}(\th) := \sup_{0<r<1}|\phi_a(re^{i\th})|$ lies in $L^p([0,2\pi])$. The conditions on $u$ imply that $u(\phi_a(re^{i\th})) \leq C|\phi_a(re^{i\th})|^{p} + C \leq C M_{\phi_a}(\th)^p+C \in L^1([0,2\pi])$. We may therefore apply the dominated convergence theorem in order to let $r \nearrow 1$ in \rrr{nowayout} and obtain

\begin{equation} \label{awayout}
u(a) \leq \lim_{r \nearrow 1} \frac{1}{2\pi} \int_{0}^{2\pi} u(\phi_a(re^{i\th}))d\th = \frac{1}{2\pi} \int_{0}^{2\pi} u(\phi_a(e^{i\th}))d\th.
\end{equation}

Recalling that $u$ is bounded above by $K$ on $\dd W$ and $\phi_a(e^{i\th}) \in \dd W$ a.e. on $[0,2\pi]$, we obtain $u(a) \leq K$. \qed

We now prove the theorem. We may assume $K=1$. Set

\begin{equation} \label{17dd}
\log^+ x = \begin{cases} \log x & x > 1,\\
0 & x  \leq 1.
\end{cases}
\end{equation}

The function $u(z) = \log^+ |f(z)|$ is the maximum of two subharmonic function,
and is therefore subharmonic. Note that the conditions on $f$ imply that
$\sup_{z \in \dd W} u(z) = 0$ and $u(z) \leq C|z|^{p} + C$ for some (possibly
different) $C>0$. Applying Proposition \ref{calig} now implies that $u(z) \leq
0$ for all $z \in W$, and the result follows. 

\section{The Green's Function of the Punctured Disk}

In the field of analysis, Green's function $G(x,y)$ on regions of $\RR^n$ is formally defined to be the solution of $LG(x,y) = \dd(y-x)$, where $L$ is a linear differential operator. The classical object in complex analysis takes $L$ to be the Laplacian, but here for a given domain $\Om \subseteq \CC$ and $z \in \Om$ the Green's function can be characterized in the following way.

\begin{defn} \label{anal}

The Green's function $G_\Om(z,w)$ on a domain $\Om$ is a function in $w$ on $\Om \bsh \{z\}$ satisfying the following properties.

\begin{itemize} \label{analdef}

\item[(i)] $G_\Om(a,z)$ is harmonic in $z$ and positive on $\Om \bsh \{a\}$.

\item[(ii)] $G_\Om(a,z) \lar 0$ as $z \lar bd(\Om)$, where $bd(\Om)$ is the set-theoretic boundary of $\Om$ (Note that the boundary is to be taken in the Riemann sphere $\hat \CC := \CC \cup \{ \ff \}$, so that if $\Om$ is unbounded then $\ff \in \Om$).

\item[(iii)] $G_\Om(a,z) + \ln |a-z|$ extends to be continuous (and therefore harmonic) at $z=a$.

\end{itemize}

\end{defn}

It can be shown by standard analytic techniques that $G_\Om$ as defined here satisfies $LG_\Om(a,z) = \dd(z-a)$. Not every domain has a Green's function as defined above, as for instance it can be shown that no function meeting these requirements can exist on the punctured disk $\DD^\times = \{0 < |z| < 1\}$, or more generally on a domain with isolated singularities. The Green's function is of tremendous importance in analysis on $\RR^n$, including complex analysis, and the question of what domains possess a Green's function has been keenly studied by analysts over the years.

\vski

On the other hand, the term ``Green's function" has entered the vocabulary of probabilists in a way that may seem initially unrelated, namely as a measure the expected number of times that a discrete process visits a point, or the expected amount of time that a continuous process spends at a point. The reason that these two different notions have garnered the same name was discovered by Hunt in 1956 (\cite{hunt}), who showed that in many cases in $\RR^n$ these two notions coincide, with $L$ the Laplacian and the process in question Brownian motion. Before discussing this fact further, let us examine the probabilistic notion of a Green's function as pertains to Brownian motion in more detail. 

\vski

We will let a Brownian motion $B_t$ start at $a$ and run until a stopping time $\tau$, and will let $B_t = \Delta$ for $t \geq \tau$, where $\Delta$ is a so-called "cemetery point" outside of $\hat \CC$; in other words, $B_t$ should no longer be in the plane for $t \geq \tau$. Let $\rho_t^\tau(a,z)$ be the probability density function at point $a$ and time $t$ of this killed Brownian motion. We then can calculate formally, for any measurable function $f$ on $\Om$,

\begin{equation} \label{exploctime}
\begin{split}
E_a \int_{0}^{\tau} f(B_s) ds & =  \int_{0}^{\tau} E_a[f(B_s)] ds \\
& = \int_{0}^{\ff} \int_{\CC} f(w) \rho_s^\tau(a,z) dA(z) ds \\
& = \int_{\CC} \Big( \int_{0}^{\ff} \rho_s^\tau(a,z) ds \Big) f(z) dA(z).
\end{split}
\end{equation}

This leads one to the consideration of the function

\begin{equation} \label{probdef}
G_\tau(a,z) := \int_{0}^{\ff} \rho_s^\tau(a,z) ds.
\end{equation}

This definition appears in a number of places, including \cite{hunt},
\cite{durBM}, \cite{bass}, and \cite{mortper}. Most often the definition has
been examined with $\tau$ being the exit time from a domain, but as is noted in
\cite{mortper} and \cite{bass} more general stopping times are allowable. Going
forward, we will always use the notation $T_\Om$ to denote the exit time of a
domain $\Om$; that is, $T_\Om = \inf\{t \geq 0 : B_t \in \Om^c\}$. In \cite{greggreenstoppingtimes}, an
expression is found for the Green's function of $T_{\DD^\times}$; note that $\DD^\times$ does not possess a Green's function in the classical sense. $G_{T_{\DD^\times}}$ was expressed as an infinite series, which was calculated using the invariance of Brownian motion under analytic functions; each term in the sum corresponded to a preimage under the covering map from the upper half-plane onto the punctured disk. It was then argued that $G_{T_{\DD^\times}} = G_{\DD}$, as with probability $1$ the puncture at the origin is avoided by Brownian motion. Equating the infinite series with the well-known Green's function of the disk resulted in a nontrivial identity, which was shown to be equivalent to Euler's infinite product identity for sine.

The purpose of this section is to show how Brownian motion can be removed from the equation; that is, to show naturally how $G_{\DD}$ can be expressed as an infinite series without resorting to
properties of Brownian motion. To begin with, it is easy to verify that the function $\ln \left|\frac{1-\ov{a}z}{a-z}\right|$ satisfies properties $(i)$-$(iii)$ above, and therefore must be equal to $G_{\DD}(a,z)$. We will prove the following.

\begin{theorem}
	The Green's function of the unit disk can be expressed as
    
	\begin{align} \label{bp}
		G_{\DD}(a,z)=\ln \left|\frac{1-\ov{a}z}{a-z}\right|=\sum_{k=-\infty}^\infty \ln\left|\frac{(\log z + 2\pi ik) + \ov{\log  a}}{(\log z + 2\pi ik) - \log a}\right| .
	\end{align}

whenever $a,z \in \DD^\times$.
\end{theorem}

A few words about the latter expression are in order. To begin with, we are distinguishing here between the logarithm of a positive real number, which we denote as $\ln$, and the analytic logarithm of a complex number, which we denote $\log$. At first glance, it seems necessary to specify which branch of $\log$ we have chosen, and perhaps to worry about a branch cut on $\DD^\times$. However, since any two branches of $\log$ differ by $2\pi i$, the choice of branch does not affect the final value of the sum, and the sum (assuming it converges) is well-defined for all $a,z$ in the punctured disk.
    
In a sense, the right side of \eqref{bp} is the Green's function of the punctured disk $\DD^\times$; at least, it is the Green's function of the exit time of Brownian motion from $\DD^\times$. Algebraically this is found by projecting the Green's
function for $\HHH$ through the covering map $w \to \exp(i w)$ and taking summation over the
preimages. However this is using a property of the Green's function which is not commonly known, namely that there should essentially be invariance under non-injective analytic functions as well as under conformal maps. In particular,
if $f:\Omega\to\Omega'$ is a covering map then

\begin{align*}
    G_{T_{\Omega'}}(f(b), f(w)) = \sum_{w'\in f^{-1}(f(w))} G_{T_\Omega}(b, w'),
\end{align*}
for any $b, w \in \Omega$; see \cite{greggreenstoppingtimes} for a probabilistic proof of this. The expression on the right side of \eqref{bp} is then obtained by applying this formula with $f(w) = \exp(i w), \Omega = \HHH$, and $\Omega' = \DD^\times$, using the fact that $G_{T_{\HHH}}(u,v) = \ln \Big|\fr{v-\ov{u}}{v-u}\Big|$ for $u,v \in \HHH$. Note that a preimage under $f$ of $a \in \DD^\times$ is given by $-i\log a$, with $\log$ any branch of the logarithm, and all preimages of $z \in \DD^\times$ are given by $-i\log z + 2 \pi k$ for $k \in \ZZ$; furthermore the conjugate of $-i\log a$ is $i \ov{\log a}$, and we have multiplied the numerator and denominator in the rightmost fraction of \eqref{bp} by $i$. 


Our goal now is to prove \eqref{bp} without any of this probabilistic machinery; that is, even though the intuition comes by way of Brownian motion, we actually don't need probability to prove it, as we now show.

\begin{proof}[Proof of Theorem 3]
In equation \eqref{bp}, it may be checked that if $z,a$ are rotated by the same amount, the value of all three expressions don't change, and so we may assume
$a$ is a positive real number. We may also chose the branch of a logarithm so that $\log a$ is a real number, making
these choices we find,
\begin{equation*}
\GG(a,z)=\sum_{k=-\infty}^\infty \ln\left|\frac{(\log z+ 2\pi ik) + \log  a}{(\log z + 2\pi ik) - \log a }\right|.
\end{equation*}
We will show that $\GG$ satisfies the conditions of the Green's function for the unit disk, which is sufficient to prove the theorem owing to the uniqueness of Green's functions. To this end, we must show that $\GG(a,z)$ extends to be a harmonic function in $z$ on $\DD$, with a logarithmic singularity as $z\to a$, and that $\GG(a,z)\to 0$ as $z$ tends to the boundary. 
We begin by showing that the series defining $\GG$ converges for fixed $a, z \in \DD^\times$. By rearranging the summation, we simplify $\GG$ to the following.

\begin{align*} 
\GG(a, z) &= \ln\left|\fr{\log z + \log  a}{\log z - \log a}\right| \\&+ \sum_{k=1}^\infty\ln\left|\left(\fr{\log z +\log  a-2\pi ik}{\log z - \log a-2\pi ik}\right)
    \left(\fr{\log z + \log  a+2\pi ik}{\log z - \log a+2\pi ik}\right)\right| \\&=
     \ln\left|\fr{\log z + \log  a}{\log z - \log a}\right| + \sum_{k=1}^\infty
    \ln\left|\fr{(\log z + \log  a)^2 +4\pi^2k^2}{(\log z - \log a)^2+4\pi^2k^2}\right|.
\end{align*} 
For simplicity, let $A(a,z) = (\log z + \log  a)^2$ and $B(a,z) = (\log z - \log a)^2$. Examining each term in the summation, we find
\begin{align}\label{eq2}
    \ln\left|\fr{A(a,z) +4\pi^2k^2}{B(a,z)+4\pi^2k^2}\right| = \ln \left|1 + \fr{A(a,z)-B(a,z)}{B(a,z) + 4\pi^2 k^2}\right|.
\end{align}
The series $\sum_{k=1}^\infty a_k$ converges absolutely if and only if $\sum_{k=1}^\infty \ln(1+a_k)$ does.
It is clear that $\sum_{k=1}^\infty\fr{A(a,z)-B(a,z)}{B(a,z) + 4\pi^2 k^2}$ converges absolutely, and so
$\GG(z,a)<\infty$.

The previous rearranging of $\GG$ also demonstrates a logarithmic singularity as $z\to w$; a singularity only occurs for the $k=0$ term, and here the nature of the singularity is clear. The quotient $z \to \frac{(\log z+ 2\pi ik) + \log  a}{(\log z + 2\pi ik) - \log a }$ is analytic in $z$, and therefore each term $z \to \ln\left|\frac{(\log z+ 2\pi ik) + \log  a}{(\log z + 2\pi ik) - \log a }\right|$ is harmonic, since following an analytic function by a harmonic one yields a harmonic composition. $\GG$ is an absolutely and locally uniformly convergent series of harmonic functions, and is therefore harmonic. Finally, it is evident that $\GG(a,z) \to 0$ as $|z| \to 1$ since a generic term in the series defining $\GG$ has the form
\begin{align} \label{kali}
    \ln\left|\frac{(\log z+ 2\pi ik) + \log  a}{(\log z + 2\pi ik) - \log a }\right|.
\end{align}
If $|z|=1$, then $\log z$ is purely imaginary, which means that the numerator and denominator in \eqref{kali} have the same imaginary parts, and real parts the same up to a sign change. As such, the modulus of the fraction is 1, and the term is equal to 0.

We now show that $\GG(a,z)$ extends to a harmonic function defined for all $z \in \DD$. Fix $a\in \DD^\times$; we need to
show that $\GG(z,a)$ is bounded in a neighborhood of the origin, then by the removable singularity theorem $\GG$ extends to a harmonic 
function on $\DD$. It is sufficient to show $\sum_{k=1}^\infty\left|\fr{A(a,z)-B(a,z)}{B(a,z) + 4\pi^2 k^2}\right|$ 
has a uniform bound as $z\to0$. For definiteness, we will let $\log z = \ln |z| + i Arg(z)$, where $Arg$ is the branch of the argument function with $-\pi < Arg(z) \leq \pi$; this is defined for all non-zero $z$, however it is clearly not continuous on the negative real axis. This is not a problem in the current situation, though, since the periodic nature of the sum defining $\GG$ results in $\GG$ being continuous, and in fact harmonic. Letting $r = \log z$, we then have,
\begin{align*}
\sum_{k=1}^\infty\left|\fr{A(a,z)-B(a,z)}{B(a,z) + 4\pi^2 k^2}\right| = \sum_{k=1}^\infty\left|\fr{4r \log a}{r^2-2r\log a + (\log a)^2 + 4\pi^2 k^2}\right|.
\end{align*}
As $z\to0$, we have that $\Re r \to -\infty$, and the imaginary part of $r$ stays in the bounded interval $(-\pi,\pi]$, so the argument of $r^2$ tends to $0$. Therefore, the magnitude of $r^2+4\pi^2 k^2$ will be close to $|r|^2+4\pi^2 k^2$, and since other terms are of lower order as $r$ become large, it is not hard to see that the series above is bounded by
\begin{align*}
    C_a\sum_{k=1}^\infty\fr{|r|}{|r|^2 + 4\pi^2 k^2} \leq C_a\int_0^{\infty} \fr{|r|}{|r|^2+4\pi^2x^2}dx = \fr{1}{4}C_a,
\end{align*}
where $C_a$ is a positive constant depending on $a$. Hence, $\GG(z,a)$ stays bounded near $z=0$, and thus extends to a harmonic function on $\DD$. We have now shown
$\GG$ satisfies all the conditions of the Green's function for the disk, and so the result follows.
\end{proof}
\section{Infinite product identities}

In this section we show how the identity we proved before yields several infinite products for trigonometric functions. Let us fix $r, c>0$ and $b \in \RR$, and then apply \eqref{bp} with $a=e^{-r}, z=e^{-c+bi}$ to get

\begin{equation} \label{}
\begin{split}
G_{\DD}(e^{-r},e^{-c+bi}) & = \ln \Big|\frac{1-e^{-r-c+bi}}{e^{-r} - e^{-c+bi}}\Big| \\
& = \sum_{n=-\ff}^{\ff} \ln \frac{|(r+c)i - (b+2\pi n)|}{|(r-c)i - (b+2\pi n)|} \\
& = \frac{1}{2} \sum_{n=-\ff}^{\ff} \ln \frac{(b+2\pi n)^2 + (r+c)^2}{(b+2\pi n)^2 + (r-c)^2}.
\end{split}
\end{equation}

Exponentiating gives the identity

\begin{equation} \label{mirror}
\prod_{n=-\ff}^\ff \frac{(b+2\pi n)^2+(r+c)^2}{(b+2\pi n)^2+(r-c)^2} = \Big|\frac{1-e^{-r-c+bi}}{e^{-r} - e^{-c+bi}}\Big|^2
\end{equation}

Special cases can take more familiar forms. If we set $b=0$ and rearrange a bit, we are led immediately to

\begin{equation} \label{}
\frac{(\frac{r+c}{2})^2}{(\frac{r-c}{2})^2} \Big(\prod_{n=1}^\ff \frac{(1+(\frac{r+c}{2\pi n})^2)}{(1+(\frac{r-c}{2\pi n})^2)}\Big)^2 = \frac{\sinh (\frac{r+c}{2})^2}{\sinh (\frac{r-c}{2})^2}.
\end{equation}

Multiplying both sides by $(\frac{r-c}{2})^2$ and taking the limit as $c \lar r$ yields the infinite product representation for $\sinh$:

\begin{equation} \label{}
\sinh r = r \prod_{n=1}^\ff (1+(\frac{r}{\pi n})^2).
\end{equation}

Returning to \rrr{mirror}, take now $b=\pi, r=c$, and we may reduce easily to

\begin{equation} \label{}
\cosh r = \prod_{n=1}^\ff (1+(\frac{r}{\pi (n-1/2)})^2),
\end{equation}

Note that the infinite product representations for sine and cosine can be derived from these, as $\sin r = i \sinh(-ir)$ and $\cos r = \cosh (ir)$, and we obtain:

\begin{equation} \label{}
\sin r = r \prod_{n=1}^\ff (1-(\frac{r}{\pi n})^2); \qquad \cos r = \prod_{n=1}^\ff (1-(\frac{r}{\pi (n-1/2)})^2).
\end{equation}

{\bf Remark:} To our knowledge, the expansion and argument given in \cite{greggreenstoppingtimes} was new, however other infinite product identities related to Green's functions can be found in \cite{melnikov2011green}.

\section{Acknowledgements} 
Part of this work was done while the first author was supported by Australian Research Council Grants DP0988483 and DE140101201.
\bibliography{CAbib}
\end{document}